\documentclass[11pt]{article}
\usepackage{amsfonts}
\usepackage{amscd}
\usepackage{amsmath}
\usepackage{epsfig}
\usepackage{color}

\usepackage{multicol}                               
\usepackage{pgf,tikz}

\usetikzlibrary{arrows, arrows.meta}

\def\be{\begin{equation}}
\def\ee{\end{equation}}
\def\ben{\begin{displaymath}}
\def\een{\end{displaymath}}
\def\baa{\begin{eqnarray}}
\def\eaa{\end{eqnarray}}

\def\ba{\begin{array}}
\def\ea{\end{array}}

\makeatletter
\@addtoreset{equation}{section}
\makeatother

\usepackage{amsthm}
\usepackage{amsmath}
\usepackage{amssymb}
\usepackage{amsbsy}
\usepackage{amsfonts}
\usepackage{latexsym}
\usepackage{amscd}
\usepackage{eucal}
\usepackage{mathrsfs}
\usepackage[all, knot]{xy}
\xyoption{arc} \xyoption{color}\xyoption{line}

\newtheorem{prop}{Proposition}
\newtheorem{lem}{Lemma}
\newtheorem{cor}{Corollary}
\newtheorem{Remark}{Remark}
\newtheorem{defi}{Definition}

\newtheorem{thm}{Theorem}

\def\be{\begin{equation}}
\def\ee{\end{equation}}
\def\ben{\begin{displaymath}}
\def\een{\end{displaymath}}
\def\baa{\begin{eqnarray}}
\def\eaa{\end{eqnarray}}

\def\ba{\begin{array}}
\def\ea{\end{array}}

\def\L{{\bf \frak L}}
\def\H{{\bf \frak H}}
\def\A{{\bf \frak A}}
\def\h{{\bf \frak h}}
\def\a{{\bf \frak a}}
\def\c{{\bf \frak c}}

\begin{document}

\title {Green function and self-adjoint Laplacians on polyhedral surfaces}

\author{ Alexey Kokotov\footnote{{\bf E-mail: alexey.kokotov@concordia.ca}}, Kelvin Lagota\footnote{{\bf E-mail: kelvin.lagota@concordia.ca}}}

\maketitle

\vskip0.5cm
\begin{center}
Department of Mathematics and Statistics, Concordia
University, 1455 de Maisonneuve Blvd. West, Montreal, Quebec, H3G
1M8 Canada \end{center}

\vskip2cm
{\bf Abstract.} Using Roelcke formula for the Green function, we explicitly construct a basis in the kernel of the adjoint Laplacian on a compact polyhedral surface $X$ and compute the $S$-matrix of $X$ at
the zero value of the spectral parameter. We apply these results to study various self-adjoint extensions of a symmetric Laplacian on a compact polyhedral surface of genus two with a single conical point.
It turns out that the behaviour of the $S$-matrix at the zero value of the spectral parameter is sensitive to the geometry of the polyhedron.

\vskip2cm

\section{Introduction}

  The spectral geometry of a Riemannian manifold  $X$ with singularities is more involved than that of smooth manifolds, in particular, due to the following reason: it may happen that the symmetric Laplacian $\Delta$ (usually defined  on smooth functions supported in $X\setminus \{singularities\}$) is not essentially self-adjoint and  in order to consider the spectrum of the Laplacian, one has to make a choice from (infinitely) many possible self-adjoint extensions of $\Delta$.

 In dimension one, this leads to the rich theory of quantum graphs; the case of Euclidean spaces ${\mathbb R}^2$ and ${\mathbb R}^3$ with punctures is investigated in great detail in \cite{Albeverio} (see also the references therein);  manifolds of higher dimension with cone like singularities are considered, e. g.,  in the papers \cite{OLDCMP}, \cite{Moors}, \cite{Kirsten}, \cite{Loya} to mention a few.
 In this paper, we consider the case of compact polyhedral surfaces (closed surfaces glued from Euclidean triangles). These are compact Riemann surfaces equipped with flat conformal metrics with conical singularities at the vertices of the corresponding polyhedron (it should be noted that the metric of a polyhedron does not see the edges:  interior points of an edge are ordinary smooth point of the corresponding Riemannian manifold).

 A question of general interest here can be formulated as follows: how do the spectral characteristics of the polyhedron depend on the choice of the self-adjoint extension of the symmetric Laplacian, the choice of conformal polyhedral metric, and moduli of the underlying Riemann surface? This question was partially addressed in \cite{HillKokJGA}, where the dependence of an important spectral invariant, the $\zeta$-regularized spectral determinant of the Laplacian,  on the choice of the self-adjoint extension was analysed. It turned out that one can write a comparison formula for two determinants of the Laplacian corresponding  to different self-adjoint extensions and the main ingredient of this formula is the so-called $S$-matrix of the polyhedral surface. The $S$-matrix depends on the spectral parameter $\lambda$ and is defined through the coefficients in the asymptotical expansions near the conical points of some special solutions (in classical sense) to the homogeneous Helmholtz equation $(\Delta-\lambda)u=0$ on the polyhedron. Moreover, the behaviour of $S(\lambda)$ at the zero value of the spectral parameter plays especially important role, say, the order of the zero of a certain minor of $S(\lambda)$ at $\lambda=0$ is related to the number of zero modes of the corresponding self-adjoint extension;  most of the entries of the matrix $S(0)$ admit explicit expression through holomorphic invariants of the underlying Riemann surface (Bergman kernel, Schiffer projective connection), and in case of smooth surface with punctures (which can be considered as conical
 points of angle $2\pi$)  some of entries of $S(0)$ are related to the Robin mass of the surface, etc.

 In the present paper, we apply and further develop the results of \cite{HillKokJGA}. In the first part of the paper, we discuss the general properties of the symmetric Laplacian $\Delta$ on arbitrary polyhedral surface: we give an explicit description of the domain of its adjoint $\Delta^*$ and, in particular, explicitly construct a basis of the kernel ${\rm Ker}\,\Delta^*$. Using the latter basis, we compute the matrix $S(0)$ expressing its entries through holomorphic invariants of the underlying Riemann surface. Our main technical tool here is the Roelcke formula for the Green function of a closed surface which we briefly discuss in the very beginning of the paper.
In the second part of the paper, we apply the results of the first part to the simplest example of a polyhedral surface  of (the lowest possible) genus two with one conical point. We study three concrete self-adjoint extensions of the symmetric Laplacian on this surface: the Friedrichs extension, the so-called holomorphic extension, and the maximal singular extension. Using the results of \cite{HillKokJGA} and the explicit formulas for $S(0)$, we write down the precise (with all the auxiliary constants computed) comparison formulas relating the $\zeta$-regularized determinants of these three extensions.  It turns out that properties of the $S$-matrix depend on geometric properties of the polyhedral surface. We show that the dimension of the kernel of the holomorphic extension (related to the order of the zero of a certain minor of $S(\lambda)$) depends on the class of linear equivalence of the divisor $(2P)$, where $P$ is the vertex of the polyhedron (this effect was previously found in \cite{HillKokAMS}, where the polyhedra of genus $g$ with $2g-2$ vertices were considered) and that the dimension of the kernel of the maximal singular extension can be higher than usual if the surface has a very large group of symmetry.

{\bf Acknowledgements.} We thank Luc Hillairet for important discussions and advices. The research of the first author was supported  by NSERC. The work was started during the stay of the first author at Max Planck Institute in Bonn, he thanks the Institute for hospitality.

\section{Green function and kernel of the adjoint Laplacian for compact polyhedral surfaces}
\subsection{Roelcke's formula for the Green function}
 Let $X$ be a compact Riemann surface and let $\rho$ be a conformal metric on $X$; we assume that $\rho$ is either smooth or flat with conical singularities. Let $\Delta^{\rho}$ be the corresponding self-adoint  Laplace operator (in the case of conical metric we define $\Delta^\rho$  as the Friedrichs extension of the symmetric Laplace operator  with domain consisting of smooth functions vanishing near the conical points: the functions from the domain of the Friedrichs extension are known to be bounded near the conical points) and let $G(x, y)$ be the corresponding Green function, i . e.,  the constant term in the expansion of the resolvent kernel, $R(x, y; \lambda)$, of  the operator $\Delta^{\rho}$ at $\lambda=0$:
 \begin{equation}\label{resolvent}
 R(x, y; \lambda)=\frac{1}{{\rm Area}\, (X)\,\lambda}+G(x, y)+O(\lambda)\,.
 \end{equation}
The Green function is real-valued and satisfies
\begin{enumerate}
 \item $G(x, y)=G(y, x)$

 \item $\Delta^{\rho}_xG(x, y)=\Delta^{\rho}_yG(x, y)=-\frac{1}{{\rm Area}(X, \rho)}$ for $x\neq y$

 \item $G(x, y)=\frac{1}{2\pi}\log|x-y|+O(1)$ as $x\to y$

 \item In the case of conical metric, the Green function $G(\, \cdot \, , y)$ is bounded near all the conical points (unless $y$ is  a conical point itself and the first argument approaches to $y$)\,.
 \item For any $x\in X$ one has \begin{equation}\label{zero}\int_X G(x, y)dS(y)=0\,,\end{equation} where $dS$ is the volume element of the metric $\rho$.
\end{enumerate}

In the case of smooth metric, the explicit formula (\ref{Roelcke}) below for the Green function is given in \cite{Fay2} (page 31, f-la (2.19)) and is called there Roelcke's formula (without any reference). Unfortunately, we were unable to identify the primary source and it seems that \cite{Fay2} is the only published text containing this result in its full generality. Formula (\ref{Roelcke}) and its proof are also valid for conical metrics. For the reader's convenience we decipher here the derivation of this formula given in passing in \cite{Fay2}.

Choosing a standard basis of $a$- and $b$-cycles on $X$ and the corresponding basis $\{v_\alpha\}_{\alpha=1}^g$ of the normalized ($\oint_{a_\beta}v_\alpha=\delta_{\alpha \beta}$) holomorphic one-forms, introduce  (see, e. g., page 4 of \cite{Fay1} with different normalization of the basic holomorphic differentials) a meromorphic
one-form $\Omega_{p-q}$ via
\begin{equation}\label{form}\Omega_{p-q}(z)=\int_p^qW(z,\, \cdot\,)-2\pi i \sum_{\alpha, \beta=1}^g\left(\Im {\mathbb B}\right)^{-1}_{\alpha \beta}v_\alpha(z)\Im \int_p^qv_\beta\,,\end{equation}
where $W$ is the canonical meromorphic bidifferential on $X$.
 It is straightforward to check that this one-form is the unique differential of the third kind with simple poles at $p$ and $q$ with residues $1$ and $-1$, and with purely imaginary periods. Thus, the real part of the integral $\int_x^y\Omega_{p-q}$ is well defined (i. e. is independent of the path of the integration) and gives a harmonic function (with logarithmic singularities) with respect to all four arguments $x, y, p, q$. Using the known singularities of the latter function, one can express it as
 \begin{equation}\label{repr}
 \Re \int_x^y\Omega_{p-q}= 2\pi\left(G(y, p) - G(y, q)+G(x, q)-G(x, p)\right).
 \end{equation}
Integrating (\ref{repr}) over $X$ twice (first with respect to $dS(x)$ and then with respect to $dS(q)$), making use of (\ref{zero}), and renaming the arguments in the resulting expression, one gets an explicit formula for the Green function
\begin{equation}\label{Roelcke}
G(x, y)=\frac{1}{2\pi ({\rm Area}\, (X, \rho))^2}\int_XdS(q)\int_XdS(p) \Re \int_p^x\Omega_{y-q}\,.
\end{equation}
\subsection{Harmonic functions with prescribed singularities}
\subsubsection{Domain of the adjoint operator and Gelfand symplectic form}
Let $P_1, \dots, P_M$  be the conical points of the metric $\rho$ and let $\beta_1, \dots, \beta_M$ be the corresponding conical angles. Introduce the integers $n_k$, \, $k=1, \dots, M$ via
$$2\pi n_k<\beta_k\leq 2\pi (n_k+1)\,.$$
 In the proof of Proposition \ref{pr2} below we will need to consider a conical point with conical angle $2\pi$: in this case $n_k=0$ and all the sums $\sum_{l=1}^{n_k}$ are equal to $0$ by definition. Let $X_0=X\setminus\{P_1, \dots, P_n\}$ and let $\Delta^*$ be the adjoint operator to the standard symmetric Laplacian $\Delta^\rho$ with domain $C^\infty_0(X_0)$. Introduce {{\it \bf the distinguished local parameter}} $\zeta_k$ near $P_k$: we remind the reader that in this local parameter one has $\rho(\zeta_k, \bar \zeta_k)|d\zeta_k|^2=
(b_k+1)^2|\zeta_k|^{2b_k}|d\zeta_k|^2$, where $2\pi(b_k+1)=\beta_k$.

In the vicinity of the point $P_k$ a function $u$ from ${\cal D}(\Delta^*)$ has the asymptotics
 \begin{equation}\label{MAINAS}u=\frac{i}{\sqrt{2\pi}}\L_k(u)\log|\zeta_k|+\sum_{m=1}^{n_k}\frac{1}{\sqrt{4\pi m}}\H_{k, m}(u)\frac{1}{\zeta_k^m}+\sum_{m=1}^{n_k}\frac{1}{\sqrt{4\pi m}}\A_{k, m}(u)\frac{1}{\bar \zeta_k^m}+\frac{i}{\sqrt{2\pi}}\c_k(u)+
 \end{equation}
 $$\sum_{m=1}^{n_k}\frac{1}{\sqrt{4\pi m}}\h_{k, m}(u)\zeta_k^m+\sum_{m=1}^{n_k}\frac{1}{\sqrt{4\pi m}}\a_{k, m}(u)\bar \zeta_k^m +\chi v,$$
 where $\chi$ is a cut-off $C^\infty$-function equal to $1$ in a small vicinity of $P_k$ with support in another small vicinity of $P_k$ and $v$ is a function from the domain of the closure ${\cal D}\overline {\Delta^\rho}$.
 One has the asymptotics $v=o(|\zeta_k|^{n_k})$ as $\zeta_k\to 0$.
 The notation for the coefficients comes from the form of the corresponding term in the asymptotics: growing holomorphic ($\H$), growing antiholomorphic ($\A$), (growing) logarithm ($\L$), constant ($\c$), and decreasing holomorphic and antiholomorphic ($\h$ and $\a$). The normalizing factors ($\frac{1}{\sqrt{4\pi m}}$, etc) are introduced to obtain the standard Darboux basis for the symplectic form in (\ref{sympl}) below.
We give the proof of (\ref{MAINAS}) in the Appendix.

 Let $\Omega$ be the symplectic form on the factor space ${\cal D}(\Delta^*)/{\cal D}({\overline \Delta^\rho})$:
 $$\Omega([u], [v]):= <\Delta^*u, \bar v>-<u, \Delta^* \bar v>,$$
 where $<u , v>=\int_Xu\bar v dS$ is the usual hermitian product. It is straightforward to show that
 \begin{equation}\label{sympl}
 \Omega([u], [v])=\sum_{k=1}^MX_k(u)\left(\begin{matrix} 0 & -I_{2n_k+1}\\ I_{2n_k+1} & 0 \end{matrix}\right)X_k(v)^t,
 \end{equation}
 where $$X(u)=$$$$(\L_k(u), \H_{k, 1}(u), \dots, \H_{k, n_k}(u), \A_{k, 1}(u), \dots, \A_{k, n_k}(u), \c_k(u), \h_{k, 1}(u), \dots, \h_{k, n_k}(u),$$ $$\a_{k, 1}(u), \dots, \a_{k, n_k}(u)).$$
\subsection{Special growing solutions and the $S$-matrix of the polyhedral surface $X$}
We remind the reader that the Friedrichs self-adjoint extension $\Delta_F$ of the operator $\Delta^\rho$ with domain $C^\infty_0(X_0)$ has the domain
\begin{equation}\label{Fried}
{\cal D}(\Delta_F)=\end{equation}
$$\{u\in {\cal D}(\Delta^*): \forall k=1, \dots, M   \ \L_k(u)=\H_{k, 1}(u)=\dots=\H_{k, n_k}=\A_{k, 1}=\dots=\A_{k, n_k}=0\}\,,
$$
in particular, all the functions from the domain of the Friedrichs extension are bounded near the conical points. The kernel of $\Delta_F$ is one dimensional and consists of constant functions.
Let $\lambda$ do not belong to the spectrum of $\Delta_F$. We define the unique {\it special growing solutions},
\begin{equation}\label{specialsol}G_{1/\zeta_k^l}(\cdot; \lambda),\ \  G_{1/\bar \zeta_k^l}(\cdot; \lambda), \ \  G_{\log|\zeta _k|}(\cdot; \lambda)\,;\end{equation}
$k=1, \dots, n; l=1, \dots, n_k$ of the homogeneous equation
\begin{equation}\label{helm}\Delta^*u-\lambda u=0\end{equation}
via their asymptotic expansions at the conical points: the unique growing term in these asymptotical expansions is that shown as the subscript of the special solution. For instance, one defines $G_{1/\zeta_k^s}$ via
$$G_{1/\zeta_k^s}(\zeta_k; \lambda)=\frac{1}{\zeta_k^s}+O(1)$$
as $\zeta_k\to 0$ and
$$G_{1/\zeta_k^s}(x; \lambda) =O(1)$$
as  $x\to P_l$ with $l\neq k$.
\begin{defi}(see \cite{HillKokJGA})
The constant terms and the coefficients near the decreasing (= positive) powers  $\zeta_l^s$ and $\bar \zeta_l^s$; $s=1, \dots, n_l$; $l=1, \dots, M$ in the asymptotic expansions of the special growing solutions form the so-called $S$-matrix, $S(\lambda)$, of the polyhedral surface $X$.
\end{defi}
Say, the entry
$S^{1/\zeta_k^r, \bar \zeta_l^s}(\lambda)$ of the $S$-matrix is given
by the coefficient near $\bar \zeta_l^s$ in the asymptotical expansion of the special growing solution $G_{1/\zeta_k^r}(\cdot; \lambda)$ near the conical point $P_l$, similarly,
the entry $S^{\log|\zeta_k|, 1_l}(\lambda)$ is  the constant term in the asymptotic expansion of the special growing solution $G_{\log |\zeta_k|}(\cdot; \lambda)$ near $P_l$.

The following proposition is a slightly improved version of Proposition 7 from \cite{HKK-Transactions}.
\begin{prop}\label{prop1} All the entries of the matrix $S(\lambda)$ except $S^{\log|\zeta_k|, 1_l}(\lambda)$ admit holomorphic continuation to $\lambda=0$, the entries  $S^{\log|\zeta_k|, 1_l}(\lambda)$ have a simple pole at $\lambda=0$.
\end{prop}
{\bf Proof.} We start with reminding the reader the construction of the special growing solutions (\ref{specialsol}).
Let $F$ be one of the following functions defined on the whole $X$
$$\chi \log|\zeta_k|,\ \ \chi \frac{1}{\zeta_k^l},\ \ \chi\frac{1}{\bar \zeta_k^l}\,,$$
where $\chi$ is a $C^\infty$ cut off function supported in a small vicinity of $P_k$ such that $\chi=1$ in some smaller vicinity of $P_k$. Let $\Delta_F$ be the Friedrichs Laplacian and let $\lambda$ do not belong to the spectrum of $\Delta_F$.  Introduce
$$f:=(\Delta^*-\lambda)F\,$$
and and define $g(\cdot; \lambda)$ as the (unique) solution of the equation
\begin{equation}\label{todiff}(\Delta_F-\lambda)g=(\Delta^*-\lambda)F\end{equation}
(it should be noticed that the right hand side of this equation belongs to $L_2(X, \rho)$).
Then
$$G(\cdot; \lambda)=F(\cdot)-g(\cdot; \lambda)$$
is the special growing solution with  principal part $F$. It follows from the above construction that
$$g(\cdot; \lambda)=g(\cdot; \lambda)+\frac{1}{{\rm Area}(X)\lambda}\int_Xf(\cdot; \lambda)-\frac{1}{{\rm Area}(X)\lambda}\int_Xf(\cdot; \lambda)=$$
\begin{equation}\label{resolv}\left[(\Delta_F-\lambda)\Big|_{1^\bot}\right]^{-1}\left((\Delta^*-\lambda)F-\frac{1}{{\rm Area}(A)}\int_X(\Delta^*-\lambda)F\right)-\frac{1}{{\rm Area}(X)\lambda}\int_Xf(\cdot; \lambda)\,.\end{equation}
The first term in (\ref{resolv}) is holomorphic in a vicinity of the point $\lambda=0$ (a simple eigenvalue of $\Delta_F$). The behaviour of the second term at $\lambda=0$ depends on the choice of the principal part $F$.
In the case $F= \chi \frac{1}{\zeta_k^l}$ or $\chi\frac{1}{\bar \zeta_k^l}$ the second term is again holomorphic at $\lambda=0$ as it is follows from the obvious relation
\begin{equation}\label{null}
\int_Xf(\cdot; 0)=0\,.
\end{equation}
If the principal part $F$ is logarithmic ($F= \chi \log|\zeta_k|$) then (\ref{null}) is no longer true
and $g(\cdot; \lambda)$ has a simple pole with residue $$-\frac{1}{{\rm Area}(X)}\int_X\Delta^*F=-\frac{2\pi}{{\rm Area}(A)}$$
Summing up, the special growing solution $G_{1/\zeta_k^l}(\cdot; \lambda)$ and $G_{1/\bar\zeta_k^l}(\cdot; \lambda)$ are holomorphic w. r. t. $\lambda$ at $\lambda=0$, whereas
 $$G_{\log|\zeta_k|}(\cdot; \lambda)=\frac{{\rm 2\pi}}{{\rm Area}(A)\lambda} + h(\cdot; \lambda)\,$$
where $h(\cdot; \lambda)$ is holomorphic near $\lambda=0$. Thus, all the coefficients in the asymptotic expansion of $G_{1/\zeta_k^l}(\cdot; \lambda)$ and $G_{1/\bar\zeta_k^l}(\cdot; \lambda)$ are holomorphic at $\lambda=0$;
the constant term in the asymptotics  $G_{\log|\zeta|}(\cdot; \lambda)$ blows up at $\lambda=0$, all other coefficients in the asymptotics  $G_{\log|\zeta|}(\cdot; \lambda)$ are holomorphic at $\lambda$=0.$\square$
\begin{Remark}The values at $\lambda=0$ of nonsingular entries  the $S$-matrix do depend on the choice of a metric $\rho$ within a given conformal class through their dependence on the distinguished local parameters of the metric
near the conical points. The opposite statement in Proposition 7 from \cite{HKK-Transactions}  was made under implicit assumption that the conformal factor is equal to one in small vicinities of the conical points.
\end{Remark}

The values of nonsingular entries of the $S$-matrix at $\lambda=0$ can be found from the asymptotics of the (unique) special growing solutions $G_{1/\zeta_k^l}(\cdot; 0)$, $G_{1/\bar \zeta_k^l}(\cdot; 0)$   of the equation
\begin{equation}\label{lapl}\Delta^*u=0\end{equation}
subject to the condition
\begin{equation}\label{ortho}\int_XudS=0\,. \end{equation}
It should be noted that there is no harmonic function on $X$ with a single logarithmic singularity, so the special growing solutions $G_{\log |\zeta_k|}(\cdot; 0)$ {\it do not exist}.
The following proposition is the first new result of the present paper.
\begin{prop}\label{pr2}

\begin{itemize}
\item The special growing solutions $G_{1/\zeta_k^l}(y; 0)$, $G_{1/\bar \zeta_k^l}(y; 0)$; $l=1, \dots, n_k$ of the equation (\ref{lapl})  are related to the coefficients of the asymptotic expansion of the Green function
 $G(\cdot, y)$ at the conical point $P_k$ via
 \begin{equation}\label{coeff}
  G(\zeta_k, y)=G(P_k, y)-\sum_{l=1}^{n_k}\frac{1}{4\pi l}G_{1/\zeta_k^l}(y; 0)\zeta_k^l-\sum_{l=1}^{n_k}\frac{1}{4\pi l}G_{1/\bar\zeta_k^l}(y; 0)\bar \zeta_k^l+o(|\zeta_k|^{n_k})\,.
 \end{equation}
\item The constant term, $G(P_k, y)$, in (\ref{coeff}) can be represented as
\begin{equation}\label{e2}
G(P_k, y)= \frac{1}{2\pi} \lim_{\lambda\to 0} \left[G_{\log|\zeta_k|}(y; \lambda)-\frac{2\pi}{{ \rm Area}(X)\lambda}\right]\,.
\end{equation}
\end{itemize}
\end{prop}
{\bf Proof.} Until the end of this proof, we assume that $X_0=X\setminus \{P_1, \dots, P_n, y\}$, i. e. we consider the point $y$ as a conical point of angle $2\pi$.
Then $G(\cdot, y)$ belongs to the domain of the operator $\Delta^*$; the latter operator is now the adjoint to the symmetric Laplacian with domain $C^\infty_0(X\setminus\{P_1, \dots, P_n, y)\}$.

 It should be noticed
that the functions $u$ from ${\cal D}(\Delta^*)$ have the asymptotics $$u(\zeta(x))=\frac{i}{\sqrt{2\pi}}\L_y(u)\log|\zeta|+\frac{i}{\sqrt{2\pi}}\c_y(u)+o(1)$$ as $x\to y$ (here the local parameter $\zeta$ is defined via $\rho=|d\zeta|^2$ near $y$ and $\zeta(y)=0$).
Since $\int_X G_{1/\zeta_k^l}dS=0$ and $\Delta^*_xG(x, y)=const=\frac{1}{{\rm Area (X)}}$, one has
\begin{equation}\label{omega}
\Omega([G(\cdot, y)], [G_{1/\zeta_k^l}(\cdot; 0)])=0\,.
\end{equation}
On the other hand, (\ref{sympl}) implies
$$ \Omega([G(\cdot, y)], [G_{1/\zeta_k^l}(\cdot; 0)])=\h_{k, l}(G(\cdot, y))\H_{k,l}(G_{1/\zeta_k^l}(\cdot; 0))-\L_y(G((\cdot, y))\c_y(G_{1/\zeta_k^l}(\cdot; 0))=$$
$$\sqrt{4\pi l}\h_{k, l}(G(\cdot; y))-\frac{1}
{\sqrt{2\pi}i}\left[ \frac{\sqrt{2\pi}}{i} G_{1/\zeta_k^l}(y; 0) \right]$$ and
$$\h_{k, l}(G(\cdot; y))=-\frac{1}{\sqrt{4\pi l}}G_{1/\zeta_k^l}(y; 0)\,.$$
Similarly,
$$\a_{k, l}(G(\cdot; 0))=-\frac{1}{\sqrt{4\pi l}}G_{1/\bar \zeta_k^l}(y; 0),$$
and (\ref{coeff}) follows.

Let $R(x, y; \lambda)$ be the resolvent kernel of $\Delta_F$. Consider the expression
\begin{equation}\label{form}
E(\lambda)=<(\Delta^*-\lambda)[R(\cdot, y; \lambda)-\frac{1}{{\rm Area}(A)\lambda}],\  G_{\log|\zeta_k|}(\cdot; \lambda)-\frac{2\pi}{{\rm Area}(X)\lambda}>-
\end{equation}
$$<R(\cdot, y; \lambda)-\frac{1}{{\rm Area}(A)\lambda}, \ (\Delta^*-\lambda)[ G_{\log|\zeta_k|}(\cdot; \lambda)-\frac{2\pi}{{\rm Area}(X)\lambda}]>
$$
Since $\lim_{\lambda\to 0}R(\cdot, y; \lambda)=G(\cdot, y)\bot 1$ and
$$\int_X [ G_{\log|\zeta_k|}(\cdot; \lambda)-\frac{2\pi}{{\rm Area}(X)\lambda}]=0$$
(the latter equality can be checked as follows
$$\int_X G=\int_X(F-g)=\int_XF-\left(\frac{1}{\lambda}\int_X\Delta^*g-\frac{1}{\lambda}\int_Xf\right)=\int_XF+\frac{1}{\lambda}\int_Xf=\frac{1}{\lambda}\int_X\Delta^*F=\frac{2\pi}{\lambda}\,),$$
one has $E(\lambda)=o(1)$ as $\lambda\to 0$.
On the other hand, computing $E(\lambda)$ via (\ref{sympl}), one gets
$$E(\lambda)=\left[G_{\log|\zeta_k|}(y; \lambda)-\frac{2\pi}{{\rm Area}(X)\lambda}\right]-2\pi\left[R(P_k, y; \lambda)-\frac{1}{{\rm Area}(X)\lambda}\right]
$$
which implies (\ref{e2}). $\square$

The next proposition immediately follows from (\ref{coeff}), \ref{e2} and Roelcke's formula (\ref{Roelcke})
\begin{prop}\label{spetialsol} One has the following explicit expressions for the special growing solutions of the homogeneous Laplace equation (\ref{lapl}) subject to (\ref{ortho})
\begin{equation}\label{otvet1}
G_{1/\zeta_k^l}(y; 0)=-\frac{1}{(l-1)! {\rm Area}\, (X)}\int_X \Omega_{y-q}^{(l-1)}(P_k)dS(q)\,\ \ \ l=1, \dots, n_k,
\end{equation}
\begin{equation}G_{1/\bar \zeta_k^l}(y; 0)=\overline {G_{1/\zeta_k^l}(y; 0)}\,.\end{equation}
Here the expression $\Omega_{y-q}^{(l-1)}(P_k)$ should be understood as follows. Write the one form $\Omega_{y-q}$ in the distinguished local parameter $\zeta_k$ in a vicinity of the conical point $P_k$:
$$\Omega_{y-q}=\omega(\zeta_k)d\zeta_k\,.$$
Then    $$\Omega_{y-q}^{(l-1)}(P_k):=\left(\frac{d}{d\zeta_k}\right)^{l-1}\omega(\zeta_k)|_{\zeta_k=0}\,.$$
Moreover, one has the relation
\begin{equation}\label{otvet2}
\lim_{\lambda\to 0}\left[G_{\log|\zeta_k|}(y; \lambda)-\frac{2\pi}{{ \rm Area}(X)\lambda}\right]=\frac{1}{{\rm Area}(X)^2}\int_X\int_X\Re \int_p^{P_k}\Omega_{y-q}\,dS(q)dS(p)\,.
\end{equation}
 \end{prop}
\subsubsection{Explicit expressions for $S(0)$}\label{sss}
Rewriting $\Omega_{y-q}$ as
\begin{equation}\label{udob}\Omega_{y-q}(z)=\int_y^qW(z,\, \cdot\,)-\pi  \sum_{\alpha, \beta=1}^g\left(\Im {\mathbb B}\right)^{-1}_{\alpha \beta}v_\alpha(z)\int_y^qv_\beta+ \pi  \sum_{\alpha, \beta=1}^g\left(\Im {\mathbb B}\right)^{-1}_{\alpha \beta}v_\alpha(z)\overline{\int_y^qv_\beta}\end{equation}
and using in (\ref{otvet2}) the reciprocity law for normalized differentials of the third kind
$$\Re\int_S^R\Omega_{P-Q}=\Re\int_Q^P\Omega_{R-S}$$
(see, e. g, \cite{FarkasKra}, p. 67),
one can easily find all the terms of the asymptotic expansions of $G_{1/\zeta_k^l}(y; 0)$ and $\lim_{\lambda\to 0} \left[G_{\log|\zeta_k|}(y; \lambda)-\frac{2\pi}{{ \rm Area}(X)\lambda}\right]$
as $y\to P_l$; $l=1, \dots, M$. This results in explicit formulas for all the finite entries of the matrix $S(0)$.
For instance, (\ref{otvet2}) and the reciprocity law immediately imply that
\begin{equation}\label{logh}
S^{\log|\zeta_k|,\  \zeta_l}(0)=\frac{1}{2{\rm Area}(X)}\int_X\Omega_{P_k-p}(P_l)dS(p);  \ l\neq k.
\end{equation}
Similarly, from (\ref{otvet1}) and (\ref{udob}), one gets the relation
\begin{equation}\label{ha}
S^{\frac{1}{\zeta_k},\  \bar \zeta_l}(0)=\pi\sum_{\alpha, \beta =1}^g(\Im {\mathbb B})^{-1}_{\alpha \beta}v_\alpha(P_k)\overline{v_\beta(P_l)}=\pi B(P_k, P_l)\,,
\end{equation}
where $B$ is the Bergman reproducing kernel for holomorphic differentials (see, e. g., \cite{Fay2}, (1.25)). (Here the the value of a differential at $P_l$ means its value in the distinguished local parameter at this point.)
Following \cite{Fay1} and \cite{Tyurin}, introduce the Schiffer bidifferential on $X$ as
\begin{equation}\label{SK}{\cal S}(P, Q)=W(P, Q)-\pi  \sum_{\alpha, \beta=1}^g\left(\Im {\mathbb B}\right)^{-1}_{\alpha \beta}v_\alpha(P)v_\beta(Q)\,.\end{equation}
The Schiffer projective connection, $S_{Sch}$, is defined via the asymptotics of the Schiffer bidifferential at the diagonal $P=Q$:
\begin{equation}\label{SCH}{\cal S}(x(P), x(Q))=\left(\frac{1}{(x(P)-x(Q))^2}+\frac{1}{6}S_{Sch}(x(P))+O(x(P)-x(Q))\right)dx(P)dx(Q)\,,\end{equation}
as $Q\to P$.
From (\ref{otvet1}) and (\ref{udob}) together with (\ref{SK}) and (\ref{SCH}), one gets (cf. \cite{HillKokAMS})
\begin{equation}\label{hh1}
S^{\frac{1}{\zeta_k},\  \zeta_l}(0)=-{\cal S}(P_k, P_l);  \ \ l\neq k,
\end{equation}
 and
 \begin{equation}\label{hh2}
 S^{\frac{1}{\zeta_k},\  \zeta_k}(0)=-\frac{1}{6}S_{Sch}(\zeta_k)\Big|_{\zeta_k=0}\,.
 \end{equation}
In the same manner, one can find explicit expressions for all the remaining (finite) entries of $S(0)$, we leave the details to the reader.
\begin{Remark}
It looks natural to define the regularized values of the singular entries of $S(\lambda)$ at $\lambda=0$ via
\begin{equation}
{\rm reg}\,S^{\log|\zeta_k|,\ 1_l}(0):=\lim_{\lambda\to 0}\left( S^{\log|\zeta_k|,\ 1_l}(\lambda)-\frac{2\pi}{{\rm Area}(X)\lambda}\right)\,.
\end{equation}
In the case of a smooth surface $X$ with a puncture $P$, considered as a conical point of angle $2\pi$ (see, e. g., \cite{CdV}, \cite{AisiouKH}),
the special growing solution $G_{\log d(\cdot, P)}(\cdot; \lambda)$ coincides with $2\pi R(\cdot, P; \lambda)$, where $R$ is the resolvent kernel of the Friedrichs extension of the Laplacian on $X\setminus P$ and $d$ is the geodesic distance on $X$; the above regularization of a (single) entry of $S(0)$,
coincides with $2\pi m(P)$, where $m(P)$ is the so-called Robin's mass (see, e. g.,  \cite{SteinerDuke}, \cite{OkCMP})
$${\rm reg}S^{\log d(P, \cdot), 1}(0)=m(P)=\lim_{Q\to P}G(P, Q)-\frac{1}{2\pi}\log d(P, Q)\,.$$
\end{Remark}

In particular, formula (\ref{otvet2}) leads to an explicit expression for $m(P)$. Unfortunately, the latter expression contains the finite part of a diverging line integral and, therefore, is not that effective as
formulas (\ref{logh}, \ref{ha}, \ref{hh1}, \ref{hh2}).
It should be noticed that using the technique of string theorists (\cite{Pol}, \cite{VerlVerl}), one can get a nice expression for the centered  Robin's mass $m(P)-\frac{M(X)}{{\rm Area}(X)}$, where $M(X)=\int_Xm(P)dS(P)$.
Following \cite{VerlVerl}, define  the function $\Phi$ on $X\times X$ via
$$-4\pi \Phi (z, w):=-2\pi\left[\int_w^z \overrightarrow{v}\right]^t(\Im {\mathbb B})^{-1}\Im  \int_w^z \overrightarrow{v} +\log\left(|E(z, w)|^2(\rho(z)\rho(w))^{1/2} \right)\,.$$
Here $\rho(z, \bar z)|dz|^2$ is the (smooth) metric on $X$ and $E(z, w)$ is the prime form (see, e. g.,  \cite{Fay1}),  $\overrightarrow{v}=(v_1, \dots, v_g)^t$. The results from \S5 of \cite{VerlVerl} imply the relation
\begin{equation}\label{VV}-G(z, w)+\frac{1}{2}m(z)+\frac{1}{2}m(w)=\Phi(z, w)+C\end{equation}
with some constant $C$.
Integrating (\ref{VV}) one gets
\begin{equation}\label{vv}\frac{M(X)}{2}+\frac{1}{2}m(w){\rm Area}(X)=\int_X\Phi(z, w)dS(z)+C{\rm Area}(X)\end{equation}
and
$$M(X){\rm Area}(X)=\int\int_{X\times X}\Phi(z, w)dS(z)dS(w)+C{\rm Area}(X)^2\,.$$
This gives the following explicit expression for  centered Robin's mass:
\begin{equation}\label{Rob}
m(w)-\frac{M(X)}{{\rm Area}(X)}=\frac{2}{{\rm Area}(X)}\int_X\Phi(z, w)dS(z)-\frac{2}{{\rm Area}(X)^2}\int\int_{X\times X}\Phi(z, w)dS(z)dS(w)\,.
\end{equation}

Moreover, from (\ref{VV}) and (\ref{vv}) follows an interesting counterpart of Roelcke formula (\ref{Roelcke})
\begin{equation}\label{GreenVV}
G(z, w)=\frac{1}{2}\left(m(z)-\frac{M(X)}{{\rm Area}(X)}\right)+\frac{1}{{\rm Area}(X)}\int_X\Phi(z, w)dS(z)-\Phi(z, w)\,
\end{equation}
mentioned in the last lines of \S5 of \cite{VerlVerl}.

\subsubsection{Kernel of $\Delta^*$} Motivated by the recent paper \cite{Shafarevich}, we will write down the basis in the kernel of the adjoint operator $\Delta^*$ (we remind the reader that $\Delta$ is the symmmetric laplacian with domain $C^\infty_0(X\setminus\{P_1, \dots, P_M\})$). This makes the constructions from Theorem 1  \cite{Shafarevich} more explicit.

Putting $v=1$ in (\ref{sympl}), one gets
\begin{equation}
\sum_{k=1}^M\frak{L}_k(u)=0
\end{equation}
for any $u\in {\rm Ker}\Delta^*$.
On the other hand, for any two points $P$ and $Q$ of $X$, there exists a harmonic function $u$ on $X\setminus\{P, Q\}$ with asymptotics
$u(x)=\log d(x, P)+O(1)$ as $x\to P$ and $u(x)=-\log d(x, Q)+O(1)$ as $x\to Q$. Thus, Proposition \ref{spetialsol} and the equality ${\rm Ker}\Delta_F=\{{\rm const}\}$ imply the following statement.
\begin{prop} The basis of ${\rm Ker}\Delta^*$ consists of
\begin{enumerate}
\item $1$
\item functions $G_{1/\zeta_k^l}( \cdot; 0)$; $k=1, \dots, M; l=1, \dots, n_k$ from Proposition \ref{spetialsol}
\item functions $G_{1/\bar \zeta_k^l}( \cdot; 0)$; $k=1, \dots, M; l=1, \dots, n_k$ from Proposition \ref{spetialsol}
\item functions $F_{P_1, P_k}(P)=\Re \int^P\Omega_{P_1-P_k}$; $k=2, \dots, M$, where $\Omega_{P_1-P_k}$ is the meromorphic one form from (\ref{form}).
\end{enumerate}
\end{prop}

\section{Self-adjoint Laplacians on genus two polyhedral surfaces with one conical point}
Here we consider several applications of the results of the previous section to concrete classes of polyhedral surfaces. In order to avoid unnecessary technical complications, we choose the simplest case of genus two surfaces with a single conical point $P$ of conical angle $6\pi$. Thus, using the setting of Section 2.2.1, one has $M=1$,  $n_1=2$, $\beta:=\beta_1=6\pi$,
\begin{equation}\label{sympl1}
 \Omega([u], [v])=X(u)\left(\begin{matrix} 0 & -I_{5}\\ I_{5} &  0 \end{matrix}\right)X(v)^t,
 \end{equation}
$$X(u)=(\L(u), \H_1(u), \H_2(u), \A_1(u), \A_2(u), \c(u), \h_1(u), \h_2(u), \a_1(u), \a_2(u)),$$
and the asymptotics in the vicinity of the point $P$ of a function $u$ from ${\cal D}(\Delta^*)$ in the distinguished local parameter $\zeta$ has the form
 $$u=\frac{1}{\sqrt{8\pi }}\H_2(u)\frac{1}{\zeta^2}
 +\frac{1}{\sqrt{8\pi }}\A_2(u)\frac{1}{\bar\zeta^2}+\frac{1}{\sqrt{4\pi }}\H_1(u)\frac{1}{\zeta}+\frac{1}{\sqrt{4\pi }}\A_1(u)\frac{1}{\bar\zeta}+\frac{i}{\sqrt{2\pi}}\L(u)\log|\zeta|
+\frac{i}{\sqrt{2\pi}}\c(u)+$$$$\frac{1}{\sqrt{4\pi }}\h_1(u)\zeta+\frac{1}{\sqrt{4\pi }}\a_1(u)\bar\zeta+
  \frac{1}{\sqrt{8\pi }}\h_2(u)\zeta^2
 +\frac{1}{\sqrt{8\pi }}\a_2(u)\bar\zeta^2
   +\chi v$$
with $v=o(|\zeta|^2)$.

We will be working with the following three self-adjoint extensions of the symmetric Laplacian $\Delta$ with domain $C^\infty_0(X\setminus\{P\})$:
\begin{itemize}
\item the Friedrichs extension, $\Delta_F$ corresponding to the lagrangian subspace of ${\rm dom}(\Delta^*)/{\rm dom}(\bar \Delta)$
$$\frak{L}(u)=\frak{H}_1(u)=\frak{H}_2(u)=\frak{A}_1(u)=\frak{A}_2(u)=0\,,$$
\item the maximal singular regular\footnote{an extension is called regular if the functions from its domain do not have logarithmic terms in the asymptotics near conical points} extension, $\Delta_{sing}$, corresponding to the lagrangian subspace
$$\frak{L}(u)=\frak{h}_1(u)=\frak{h}_2(u)=\frak{a}_1(u)=\frak{a}_2(u)=0\,,$$
\item the holomorphic extension, $\Delta_{hol}$, corresponding to the lagrangian subspace
 $$\frak{L}(u)=\frak{A}_1(u)=\frak{A}_2(u)=\frak{a}_1(u)=\frak{a}_2(u)=0\,.$$
\end{itemize}
\begin{prop} The operators $(\Delta_{sing}-\lambda)^{-1}-(\Delta_F-\lambda)^{-1}$ and $(\Delta_{hol}-\lambda)^{-1}-(\Delta_F-\lambda)^{-1}$ are finite dimensional and one has the following representations for their traces:
\begin{equation}\label{T}
{\rm Trace}\,[(\Delta_{sing}-\lambda)^{-1}-(\Delta_F-\lambda)^{-1}]=
-{\rm Trace}\left( T^{-1}(\lambda)T'(\lambda)\right)\,
\end{equation}
\begin{equation}\label{P}
{\rm Trace}\,[(\Delta_{hol}-\lambda)^{-1}-(\Delta_F-\lambda)^{-1}]=-{\rm Trace}\left( P^{-1}(\lambda)P'(\lambda)\right)\,,
\end{equation}
where the matrices $T(\lambda)$ and $P(\lambda)$ are given in (\ref{matT}) and (\ref{matP}) below.
\end{prop}
{\bf Proof.}
Notice that the kernel of the operator $\Delta^*-\lambda$ with $\lambda\in {\mathbb C}\setminus {\rm Spectrum}(\Delta_F)$ is generated by the special growing solutions
$$G_{1/\zeta^2}(\cdot; \lambda), \ G_{1/\bar \zeta^2}(\cdot; \lambda), \ G_{1/\zeta}(\cdot; \lambda), \ G_{1/\bar \zeta}(\cdot; \lambda), \ G_{\log|\zeta|}(\cdot; \lambda)$$
 of the equation $\Delta^*u-\lambda u=0$
and, therefore, the deficiency indices of $\Delta$ are $(5, 5)$. So, we are in a position to use the Krein formula for the difference of the resolvents of two self-adjoint extensions of a symmetric operator with (equal) finite deficiency indices (see, e. g.,  \S 84 of \cite{AkhGlaz}):
\begin{equation}\label{Krein}
[(\Delta_{sing}-\lambda)^{-1}-(\Delta_F-\lambda)^{-1}](f)=\end{equation}$$\sum_{\alpha=1/\zeta^2,\ 1/\zeta, \ 1/\bar \zeta^2,\ 1/\bar \zeta}G_\alpha(\cdot; \lambda)
\sum_{\beta= 1/\zeta^2,\ 1/\zeta,\ 1/\bar \zeta^2,\ 1/\bar \zeta, \ \log|\zeta|} x_{\alpha \beta}(\lambda)<f, G_\beta (\cdot, \bar \lambda)>=$$
$$\sum_{\alpha= 1/\zeta^2,\ 1/\zeta,\ 1/\bar \zeta^2,\ 1/\bar \zeta} X_{\alpha}(\lambda)G_\alpha(\cdot; \lambda)
$$
Introducing $u\in {\rm dom\,} (\Delta_F)$ via $(\Delta_F-\lambda)u=f$ and comparing the coefficients in the asymptotic expansion of the left and right hand sides of (\ref{Krein}), one gets
\begin{equation}\label{linsys}(\frac{1}{\sqrt{4\pi}}\frak{h}_1(u), \frac{1}{\sqrt{8\pi}}\frak{h}_2(u),\frac{1}{\sqrt{4\pi}}\frak{a}_1(u),\frac{1}{\sqrt{8\pi}}\frak{a}_2(u))^t=T(\lambda)(X_{1/\zeta}, X_{1/\zeta^2}, X_{1/\bar \zeta}, X_{1/\bar \zeta^2})^t\end{equation}
with
\begin{equation}\label{matT}T(\lambda)=\left(\begin{matrix}
S^{1/\zeta,\  \zeta}(\lambda)\ \ S^{1/\zeta^2, \  \zeta}(\lambda) \ \ S^{1/\bar \zeta,\  \zeta}(\lambda)\ \ S^{1/\bar \zeta^2,\  \zeta}(\lambda) \\
S^{1/\zeta,\  \zeta^2}(\lambda)\ \ S^{1/\zeta^2,\   \zeta^2}(\lambda) \ \ S^{1/\bar \zeta,\  \zeta^2}(\lambda)\ \ S^{1/\bar \zeta^2,\  \zeta^2}(\lambda) \\
S^{1/\zeta, \ \bar \zeta}(\lambda)\ \ S^{1/\zeta^2,\   \bar \zeta}(\lambda) \ \ S^{1/\bar \zeta,\  \bar \zeta}(\lambda)\ \ S^{1/\bar \zeta^2,\  \bar \zeta}(\lambda) \\
S^{1/\zeta,\  \bar \zeta^2}(\lambda)\ \ S^{1/\zeta^2, \  \bar \zeta^2}(\lambda) \ \ S^{1/\bar \zeta,\  \bar \zeta^2}(\lambda)\ \ S^{1/\bar \zeta^2,\  \bar \zeta^2}(\lambda)
\end{matrix}
\right)
\end{equation}
Since (\ref{linsys}) holds with an arbitrary left hand side (one can take as $u$ an arbitrary function from ${\rm dom}(\Delta_F)$), the matrix $T(\lambda)$ is invertible.

Notice that
$$<(\Delta^*-\lambda)u, G_{1/\bar \zeta}(\cdot; \bar \lambda)>=<(\Delta^*-\lambda)u, \overline{G_{1/\zeta}(\cdot; \lambda)}>-<u, \overline{(\Delta^*-\lambda)G_{1/\zeta}(\cdot; \lambda)}>$$
$$=\Omega(u, G_{1/ \zeta}(\cdot; \lambda))=\sqrt{4\pi}\frak{h}_1(u)\,.$$
Similarly,
\begin{equation}\label{coeffform}
<(\Delta^*-\lambda)u, G_{1/\zeta}(\cdot; \bar \lambda)>=\sqrt{4\pi}\frak{a}_1(u); \end{equation} $$ <(\Delta^*-\lambda)u, G_{1/\bar\zeta^2}(\cdot; \bar \lambda)>=\sqrt{8\pi}\frak{h}_2(u)$$
and
$$<(\Delta^*-\lambda)u, G_{1/\zeta^2}(\cdot; \bar \lambda)>=\sqrt{8\pi}\frak{a}_2(u)\,.$$
Differentiating (\ref{todiff}) with respect to $\lambda$ and using (\ref{coeffform}), one gets the relations
$$\frac{d}{d\lambda}S^{*, \zeta}(\lambda)=-\frac{1}{4\pi}<G_{*}(\cdot, \lambda), G_{1/\bar \zeta}(\cdot; \bar \lambda)>\,,$$
$$\frac{d}{d\lambda}S^{*, \zeta^2}(\lambda)=-\frac{1}{8\pi}<G_{*}(\cdot, \lambda), G_{1/\bar \zeta^2}(\cdot; \bar \lambda)>\,,$$
and
$$\frac{d}{d\lambda}S^{*, \bar \zeta}(\lambda)=-\frac{1}{4\pi}<G_{*}(\cdot, \lambda), G_{1/ \zeta}(\cdot; \bar \lambda)>\,,$$
\begin{equation}\label{Sdiff}\frac{d}{d\lambda}S^{*, \bar \zeta^2}(\lambda)=-\frac{1}{8\pi}<G_{*}(\cdot, \lambda), G_{1/\zeta^2}(\cdot; \bar \lambda)>\,,\end{equation}
 with $$*=1/\zeta, 1/\bar \zeta, 1/\zeta^2, 1/\bar\zeta^2\,.$$
Now (\ref{Krein}) can be rewritten as
\begin{equation}\label{fortrace}[(\Delta_{sing}-\lambda)^{-1}-(\Delta_F-\lambda)^{-1}](f)=\end{equation}
$$(G_{1/\zeta}(\cdot; \lambda), G_{1/\zeta^2}(\cdot; \lambda), G_{1/\bar\zeta}(\cdot; \lambda), G_{1/\bar \zeta^2}(\cdot; \lambda))\times $$$$T^{-1}(\lambda)(\frac{1}{\sqrt{4\pi}}\frak{h}_1(u), \frac{1}{\sqrt{8\pi}}\frak{h}_2(u),\frac{1}{\sqrt{4\pi}}\frak{a}_1(u),\frac{1}{\sqrt{8\pi}}\frak{a}_2(u))^t=$$
$$(G_{1/\zeta}(\cdot; \lambda), G_{1/\zeta^2}(\cdot; \lambda), G_{1/\bar\zeta}(\cdot; \lambda), G_{1/\bar \zeta^2}(\cdot; \lambda))T^{-1}(\lambda)\times $$
$$
(\frac{1}{4\pi}<f, G_{1/\bar \zeta}(\cdot, \bar \lambda)>, \frac{1}{8\pi}<f, G_{1/\bar\zeta^2}(\cdot, \bar\lambda)>,\frac{1}{4\pi}<f, G_{1/ \zeta}(\cdot; \bar \lambda)>, \frac{1}{8\pi}<f, G_{1/\zeta^2}(\cdot; \bar\lambda)>)^t\,.$$
Relation (\ref{T}) immediately follows from (\ref{fortrace}), the elementary relation $${\rm Trace}\, g<\cdot, h>=<g, h>$$ and identities (\ref{Sdiff}).

Similarly,
\begin{equation}\label{Krein-hol}
[(\Delta_{hol}-\lambda)^{-1}-(\Delta_F-\lambda)^{-1}](f)=\sum_{\alpha= 1/\zeta^2,\ 1/\zeta} X_{\alpha}(\lambda)G_\alpha(\cdot; \lambda)
\end{equation}
and
\begin{equation}\label{linsys-2}(\frac{1}{\sqrt{4\pi}}\frak{a}_1(u), \frac{1}{\sqrt{8\pi}}\frak{a}_2(u))^t=P(\lambda)(X_{1/\zeta}(\lambda), X_{1/\zeta^2}(\lambda))^t\end{equation}
with
\begin{equation}\label{matP}P(\lambda)=\left(\begin{matrix}
S^{1/\zeta, \ \bar \zeta}(\lambda)\ \ S^{1/\zeta^2,\   \bar \zeta}(\lambda) \\
S^{1/\zeta,\  \bar \zeta^2}(\lambda)\ \ S^{1/\zeta^2, \  \bar \zeta^2}(\lambda)
\end{matrix}
\right)
\end{equation}
and (\ref{P}) follows from the same considerations as above. $\square$

The following Proposition is an immediate corollary of (\ref{otvet1}) (cf. Section \ref{sss}).
\begin{prop}\label{szero}
Introduce the function $H(\cdot, \cdot)$ (both arguments are distinguished local parameters in a small vicinity of $P$) via
$$W=\left[\frac{1}{(\zeta(Q)-\zeta(R))^2} +H(\zeta(Q),\zeta(R)\right]d\zeta(Q)d\zeta(R)$$ as $Q, R \to P$,  where $W$ is the canonical meromorphic bidifferential on $X$
(in particular, one has the relation
$$6H(\zeta(P), \zeta(P))=S_B(\zeta(P))\,,$$
where $S_B$ is the Bergman projective connection).
Then the matrix $T(0)$ is given via
\begin{equation}\label{TT}
\left(\begin{matrix}T_{11}(0) & T_{12}(0)\\
T_{21}(0) & T_{22}(0)\\
T_{31}(0) & T_{32}(0)\\
T_{41}(0) & T_{42}(0)
   \end{matrix}\right)=\end{equation} $$\left(\begin{matrix}-\frac{1}{6}S_{Sch}(0) & -H'_{\zeta}(\zeta, \zeta')|_{(0,0)}+\pi\sum(\Im {\mathbb B})^{-1}_{\alpha \beta}v_\alpha(0)v'_\beta(0)\\
-\frac{1}{2}H'_{\zeta'}(\zeta, \zeta')|_{(0, 0)}+\frac{\pi}{2}\sum(\Im {\mathbb B})^{-1}_{\alpha \beta}v_\alpha(0)v'_\beta(0) & -\frac{1}{2}H''_{\zeta \zeta'}(\zeta, \zeta')|_{(0, 0)}+\frac{\pi}{2}\sum(\Im {\mathbb B})^{-1}_{\alpha \beta}v'_\alpha(0)v'_\beta(0)   \\
\pi B(0, 0) & \pi \sum(\Im {\mathbb B})^{-1}_{\alpha \beta}v'_\alpha(0)\overline{v_\beta(0)} \\
\frac{\pi}{2} \sum(\Im {\mathbb B})^{-1}_{\alpha \beta}v_\alpha(0)\overline{v'_\beta(0)} & \frac{\pi}{2}\sum(\Im {\mathbb B})^{-1}_{\alpha \beta}v'_\alpha(0)\overline{v'_\beta(0)}
   \end{matrix}\right)
$$
and
\begin{equation}
\left(\begin{matrix}T_{13}(0) & T_{14}(0)\\
T_{23}(0) & T_{24}(0)\\
T_{33}(0) & T_{34}(0)\\
T_{43}(0) & T_{44}(0)
   \end{matrix}\right)=
   \left(\begin{matrix}\overline{ T_{31}(0)} & \overline{ T_{32}(0)}\\
\overline{ T_{41}(0)} & \overline{ T_{42}(0)}\\
\overline{T_{11}(0)} & \overline{ T_{12}(0)}\\
\overline{ T_{21}(0)} & \overline{ T_{22}(0)}
                                                                    \end{matrix}\right).
\end{equation}

One also has
\begin{equation}\label{PP}
P(0)=\left(\begin{matrix}
T_{31}(0) & T_{32}(0)\\
T_{41}(0) & T_{42}(0)
\end{matrix}\right).
\end{equation}

\end{prop}
The next proposition describes the asymptotic behaviour of the $S$-matrix as $\lambda \to -\infty$.
\begin{prop} All the entries of the matrix $T(\lambda)$ except $S^{1/\zeta, \bar \zeta}(\lambda)$, $S^{1/\zeta^2, \bar \zeta^2}(\lambda)$ and their conjugate  $S^{1/\bar \zeta,  \zeta}(\lambda)$, $S^{1/\bar \zeta^2,  \zeta^2}(\lambda)$ are $O(|\lambda|^{-\infty})$ as $\lambda\to -\infty$. One has the asymptotics
\begin{equation}S^{1/\zeta, \bar \zeta}(\lambda)=-\frac{2^{1/3}\sqrt{3}\Gamma(2/3)}{\pi\Gamma(4/3)}(-\lambda)^{1/3}+O(|\lambda|^{-\infty});
\end{equation}
$$S^{1/\zeta^2, \bar \zeta^2}(\lambda)=-\frac{2^{-1/3}\sqrt{3}\Gamma(1/3)}{\pi\Gamma(5/3)}(-\lambda)^{2/3}+O(|\lambda|^{-\infty})
$$
and
\begin{equation}\label{asdet}
{\rm det}\,T(\lambda) = \left(\frac{27}{2\pi^2}\right)^2 \lambda^2+O(|\lambda|^{-\infty})          \ \ \ \ \ \  {\rm det}\,P(\lambda)=-\frac{27}{2\pi^2} \lambda+O(|\lambda|^{-\infty})
\end{equation}
as $\lambda\to -\infty$.
\end{prop}
{\bf Proof.} (cf. \cite{HKK-Transactions}.) Passing to polar coordinates, $r, \phi$ such that $\zeta=r^{1/3}e^{i\phi/3}$; $0\leq \phi \leq 6\pi$,
one finds that the functions
$$K_\nu(\sqrt{-\lambda}r)e^{-i\nu\phi}; \ \ \ \nu= \frac{1}{3},\frac{2}{3}
 \,,$$
 where $K_\nu$ is the modified Bessel function, satisfy
 the equation (\ref{helm})
 in a vicinity of $P$.
  The well-known asymptotics of the modified Bessel function (with  $\nu>0$) reads as
  $$K_\nu(y)=\frac{\pi}{2\sin(\nu \pi)}\left[\frac{y^{-\nu}}{2^{-\nu}\Gamma(1-\nu)}-  \frac{y^\nu}{2^\nu\Gamma(1+\nu)}+O(y^{2-\nu})\right]$$
  as $y\to 0$.
  Thus, the functions $\Phi_\nu:=\pi^{-1}2^{-\nu}\Gamma(1-\nu)\sin(\pi\nu)(\sqrt{-\lambda})^\nu K_\nu(\sqrt{-\lambda}r)e^{-i\nu\phi}$; $\nu=1/3, 2/3$
  satisfy (\ref{helm}) in a vicinity of $P$ and have the asymptotics
  \begin{equation}\label{Phi1}
  \Phi_{1/3}(\zeta, \bar \zeta; \lambda)=\frac{1}{\zeta}-\frac{2^{1/3}\sqrt{3}\Gamma(2/3)}{\pi\Gamma(4/3)}(-\lambda)^{1/3}\bar\zeta +o(|\zeta|^2)
  \end{equation}
  $$\Phi_{2/3}(\zeta, \bar \zeta; \lambda)=\frac{1}{\zeta^2}-\frac{2^{-1/3}\sqrt{3}\Gamma(1/3)}{\pi\Gamma(5/3)}(-\lambda)^{2/3}\bar\zeta^2 +o(|\zeta|^2)
  $$
as $\zeta\to 0$.

Now, notice that one can change the construction of the special growing solutions from the proof of  Proposition \ref{prop1} replacing function $F$ by $\Phi_\nu$; this gives
\begin{equation}\label{altspec}G_{1/\zeta}(\cdot; \lambda)=\Phi_{1/3}(\cdot; \lambda)-(\Delta_F-\lambda)^{-1}(\Delta^*-\lambda)[\chi \Phi_{1/3}(\cdot; \lambda)]; \end{equation}

$$G_{1/\zeta^2}(\cdot; \lambda)=\Phi_{2/3}(\cdot; \lambda)-(\Delta_F-\lambda)^{-1}(\Delta^*-\lambda)[\chi \Phi_{2/3}(\cdot; \lambda)]\,.$$

Since $K_\nu(x)$ and all its derivatives are $O(e^{-x})$ as $x\to +\infty$ and
 the support of $(\Delta^*-\lambda)[\chi \Phi_{\nu}(\cdot; \lambda)]$ is separated from the origin, all the coefficients in the asymptotic expansions (\ref{MAINAS}) of second terms in the right hand sides of  (\ref{altspec}) are exponentially decreasing as $\lambda\to -\infty$ and, therefore, all the statements of the proposition follow from (\ref{Phi1}).

The next proposition is a direct consequence of Theorem 2 from \cite{HillKokJGA} and (\ref{asdet}).
\begin{prop}
Introduce the zeta-regularized determinants of the operators $\Delta_F-\lambda$, $\Delta_{sing}-\lambda$ and $\Delta_{hol}-\lambda$ via
$${\rm det}\, A= \exp\{- \zeta'_A(0)\},$$
 where $\zeta_A(s)$ is the operator zeta-function of an operator $A$ (without zero modes). Then
 \begin{equation}\label{F-sing}
 {\rm det}(\Delta_{sing}-\lambda)= \left(\frac{2\pi^2}{27}\right)^2 {\rm det}\,T(\lambda){\rm det}(\Delta_{F}-\lambda)
 \end{equation}
for real $\lambda$ not belonging to the union of the spectra of $\Delta_F$ and $\Delta_{sing}$. Similarly
 \begin{equation}\label{F-hol}
 {\rm det}(\Delta_{hol}-\lambda)= \frac{2\pi^2}{27} {\rm det}\,P(\lambda){\rm det}(\Delta_{F}-\lambda)
 \end{equation}
for real $\lambda$ not belonging to the union of the spectra of the operators
$\Delta_F$ and $\Delta_{hol}$.
\end{prop}

Since ${\rm dim}\, {\rm Ker}\, \Delta_F=1$,  the above Proposition shows that the order of the zero of ${\rm det}\, T(\lambda)$   (resp. ${\rm det}\, P(\lambda)$)  at $\lambda=0$ is one unit less than the dimension of the kernel of $\Delta_{sing}$ (resp. $\Delta_{hol}$).
We will prove in the next subsection that generically ${\rm dim}\,{\rm Ker} \Delta_{hol}=1$. We conjecture that this is also the case for $\Delta_{sing}$ (i. e. that generically ${\rm det} T(0)\neq 0$). However, we will show that by choosing  a ``very symmetric'' polyhedron $X$, one can get ${\rm dim}\,{\rm Ker} \Delta_{sing}=3$.

So, under assumption of genericity, passing to the limit $\lambda\to 0$ in (\ref{F-sing}) and (\ref{F-hol}), one gets the following comparison formulas for modified (i. e. with zero modes excluded) determinants of s. a. extensions $\Delta_F$, $\Delta_{sing}$ and $\Delta_{hol}$.
\begin{thm} Let ${\rm dim}\,{\rm Ker}\,\Delta_{sing}=1$.
Then
 \begin{equation}\label{F-sing-det}
 {\rm det}^*\Delta_{sing}= \left(\frac{2\pi^2}{27}\right)^2 {\rm det}\,T(0){\rm det}^*\Delta_{F},
 \end{equation}
 where $T(0)$ is explicitly given in (\ref{TT}).

 Let ${\rm dim}\,{\rm Ker}\,\Delta_{hol}=1$ (i. e. $P$ is not a Weierstrass point of $X$, see Proposition \ref{Weier} below). Then
 \begin{equation}\label{F-hol-det}
 {\rm det}^*\Delta_{hol}= \frac{2\pi^2}{27} {\rm det}\,P(0){\rm det}^*\Delta_{F},
 \end{equation}
 where $P(0)$ is given in (\ref{PP}).
\end{thm}
\begin{Remark} If $(2P)=C$, where $C$ is the canonical class, the flat metric on $X$ with a single conical point at $P$ has the form $|\omega|^2$, where $\omega$ is a holomorphic differential on $X$ with double zero at $P$.
In this case, an explicit expression for ${\rm det}^*\Delta_{F}$ can be found in \cite{KK-JDG}. An explicit formula for   ${\rm det}^*\Delta_{F}$  for an arbitrary $P$ can be found in \cite{Kok-PAMS}.
\end{Remark}

\begin{Remark}Let us mention two geometric constructions leading to a flat surface $X$ of genus two with a single conical singularity.
\begin{enumerate}
{\rm
\item Take a compact Riemann surface $X$  of genus two and choose a point $P\in X$. Then according to the  Troyanov theorem (\cite{Troyanov}), there exists the unique (up to rescaling) flat conformal metric on $X$ with conical singularity of angle $6\pi$ at $P$. If the divisor $(2P)$ is in the canonical class, then there exists a holomorphic one form on $X$ with divisor $2P$ and the Troyanov metric necessarily coincides (up to resacaling) with
    $|\omega|^2$. In this case, the metric has trivial holonomy. If the divisor $(2P)$ does not belong to the canonical class, then the Troyanov metric must have nontrivial holonomy along some nontrivial cycle on $X$.
    (It should be noted that the holonomy of the Troyanov metric along a small loop around the conical point is always trivial: the tangent vector turns to the angle $6\pi$ after parallel transform along this loop.)

\item (See Figure 1.) In case of trivial holonomy, the flat surface $X$ can be produced via the well known pentagon construction (see, e. g., \cite{McMullen}). Consider a pentagon $\Pi$ in the complex plane. Let the center of one of its sides coincide with the origin. Gluing the parallel sides of the octagon $\Pi\cup(-\Pi)$ together one gets a flat surface $X$ of genus $2$ with a single conical singularity of conical angle $6\pi$. The one form $dz$ in the complex
plane gives rise to a holomorphic one form $\omega$ on $X$ with a single double zero at the point $P$ on $X$ that came from eight vertices of the octagon glued together. The natural flat metric on $X$ has trivial holonomy and coincides with $|\omega|^2$.

Now take the octagon $\Pi\cup(-\Pi)$ and deform it keeping the lengths of all the sides fixed (after this deformation the opposite sides are no longer parallel). Glue the sides together following the same gluing scheme as before. Again one gets a flat surface of genus two with a single conical singularity of angle $6\pi$ but now the corresponding flat metric has nontrivial holonomy: the parallel transport along  the closed loop which came from a segment connecting two points on the opposite sides of the deformed octagon turns the tangent vector for the angle which is equal to the angle between these two opposite sides.}

\end{enumerate}

\end{Remark}

\begin{figure}[h]
\begin{center}
\begin{multicols}{2}


\begin{tikzpicture}[scale=0.5,line cap=round,line join=round,>=triangle 45,x=1.0cm,y=1.0cm]

\draw[-,color=black] (-2,0) -- (2,1);
\draw[color=black] (0,.5) node[below] {$d$};
\draw[-,color=black] (-2,0) -- (-3,3);
\draw[color=black] (-2.5,1.5) node[left] {$c$};
\draw[-,color=black] (-3,3) -- (-1,6);
\draw[color=black] (-2,4.5) node[left] {$b$};
\draw[-,color=black] (-1,6) -- (2,4);
\draw[color=black] (0.5,5) node[above] {$a$};
\draw[-,color=black] (2,4) -- (2,1);
\draw[color=black] (2,2.5) circle (1.5pt);
\draw[color=black] (2,2.5) node[right] {\tiny $O$};

\draw[-,color=black] (5,-1) -- (2,1);
\draw[color=black] (3.5,0) node[below] {$a$};
\draw[-,color=black] (5,-1) -- (7,2);
\draw[color=black] (6,0.5) node[right] {$b$};
\draw[-,color=black] (7,2) -- (6,5);
\draw[color=black] (6.5,3.5) node[right] {$c$};
\draw[-,color=black] (2,4) -- (6,5);
\draw[color=black] (4,4.5) node[above] {$d$};

\end{tikzpicture}


\begin{tikzpicture}[scale=0.5,line cap=round,line join=round,>=triangle 45,x=1.0cm,y=1.0cm]

\draw[-,color=black] (2,5) -- (-1,7);
\draw[color=black] (0.5,6) node[above] {$a$};
\draw[-,color=black] (-4,5) -- (-1,7);
\draw[color=black] (-2.5,6) node[above] {$b$};
\draw[-,color=black] (-4,5) -- (-2.103796,2.46931);
\draw[color=black] (-3.1,3.75) node[left] {$c$};
\draw[-,color=black] (-2.103796,2.46931) -- (2,1);
\draw[color=black] (0,1.75) node[below] {$d$};

\draw[-,color=black] (4,-2) -- (2,1);
\draw[color=black] (3,-0.5) node[left] {$a$};
\draw[-,color=black] (4,-2) -- (6,1);
\draw[color=black] (5,-0.5) node[right] {$b$};
\draw[-,color=black] (6.03706,4.16206) -- (6,1);
\draw[color=black] (6,2.5) node[right] {$c$};
\draw[-,color=black] (6.03706,4.16206) -- (2,5);
\draw[color=black] (4,4.5) node[above] {$d$};

\end{tikzpicture}

\end{multicols}
\end{center}

\caption{Gluing schemes for $X$: trivial (left) and nontrivial (right) holonomy}
\end{figure}
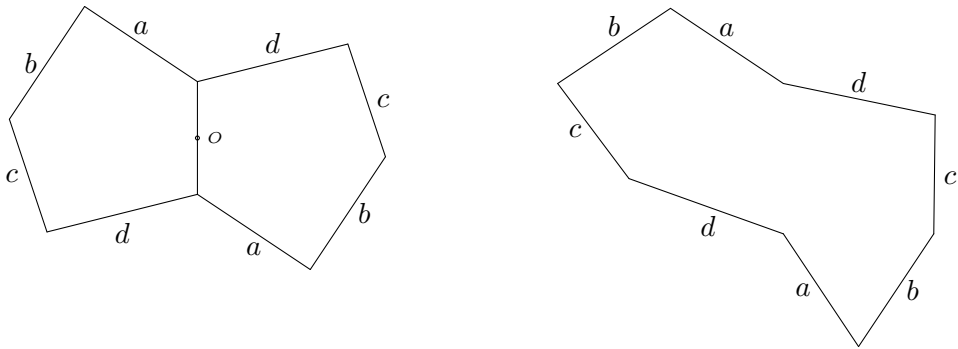

\subsubsection{One more comparison formula for resolvent kernels.} Here we briefly describe an interesting counterpart to formula (\ref{Krein-hol}) which holds in case of general conformal flat conical metrics of trivial holonomy on compact Riemann surfaces $X$ of genus $g\geq 2$. All these metrics have the form $|\omega|^2$, where $\omega$ is a holomorphic one form on $X$. Flat surfaces $X$ of genus $2$ with a single conical point $P$ of angle $6\pi$ enter this class if and only if $P$ is a Weierstrass point of $P$.

\begin{prop} Let the metric on $X$ be given by $|\omega|^2$, where $\omega$ is a holomorphic one form. Let $P_1, P_2, \dots, P_{M};    M\leq 2g-2$ be the distinct  zeroes of $\omega$ or, what is the same, the conical points of the metric $|\omega|^2$.
 Then there is the following relation between the resolvents, $R_{hol}$ and $R_F$, of the holomorphic and Friedrichs  extensions of the symmetric Laplacian on $X\setminus\{P_1, \dots, P_{M}\}$:
\begin{equation}\label{comp2}
R_{hol}(x, y; \lambda)=\frac{4}{\lambda} \frac{1}{\omega(x)\overline{\omega(y)}}\partial_x\partial_{\bar y}R_F(x, y; \lambda)
\end{equation}
\end{prop}

{\bf Proof.} We start with reminding the reader the standard relation
\begin{equation}\label{Bkernel}
\partial_x\partial_{\bar y}G_F(x, y)=-\frac{1}{4}\sum_{\alpha, \beta=1}^g\left(\Im {\mathbb B}  \right)_{\alpha \beta}^{-1}v_\alpha(x)\overline{v_\beta(y)}=-\frac{1}{4}B(x, \bar y),
\end{equation}
where $B(x, \bar y)$ is the reproducing kernel for holomorphic differentials. Here $G_F$ is just the Green function from (\ref{Roelcke}), the subscript is introduced to emphasize that we deal with the Green function of the Friedrichs Laplacian.
 Equation (\ref{Bkernel}) directly follows from (\ref{Roelcke}) (the factor $1/4$ appears due to the presence of the factor $4$ in the definition of the Laplacian, some authors do not introduce these factors).

According to  \cite{Hillairet}, one has the relations
\begin{equation}\label{DtN}
\Delta_F=4 D_z^*D_z, \ \ \ \ \ \Delta_{hol}=4D_zD_z^*;
\end{equation}
where $D_z$ is the closure of the operator
$$\frac{1}{\omega}\partial_z: C^\infty_0(X\setminus \{P_1, \dots, P_{M}\})\subset L_2(X, |\omega|^2)\longrightarrow L_2(X, |\omega|^2)\,.$$
Clearly, $D^*_z$ acts as $\frac{1}{\bar \omega}\partial_{\bar z}$.

 Now (\ref{DtN}) immediately implies that the function $\phi_m$ is a normalized eigenfunction of $\Delta_F$ corresponding to a nonzero eigenvalue $\lambda_m$ if and only if $\frac{2}{\sqrt{\lambda_m}}D_z\phi_m$ is a normalized eigenfunction of $\Delta_{hol}$ corresponding the eigenvalue $\lambda_m$.
Taking into account that  ${\rm Ker}\,\Delta_{hol}$ is spanned by the functions $\frac{v_\alpha}{\omega}$ and, therefore, the orthogonal projection in $L_2(X, |\omega|^2)$ onto ${\rm Ker}\,\Delta_{hol}$ is the integral operator with the integral kernel $\frac{B(x, \bar y)}{\omega(x)\overline{\omega(y)}}$, one gets
the following representation for the resolvent kernel of $\Delta_{hol}$ (in the sense of distribution theory):
\begin{equation}
R_{hol}(x, y; \lambda)=-\frac{B(x, \bar y)}{\omega(x)\overline{\omega(y)}}\frac{1}{\lambda}+4\sum_{\lambda_m\neq 0}\frac{D_x\phi_m(x)\overline{D_y\phi_m(y)}}{(\lambda_m-\lambda)\lambda_m}\,.
\end{equation}
Taking into account the relations $\frac{1}{(\lambda_m-\lambda)\lambda_m}=\frac{1}{\lambda}\left(\frac{1}{\lambda_m-\lambda}-\frac{1}{\lambda_m}\right)\,,$
$R_F(x, y, \lambda)=-\frac{1}{A\lambda}+\sum_{\lambda_m\neq 0}\frac{\phi_m(x)\overline{\phi_m(y)}}{\lambda_m-\lambda}\,,$
$G_F(x, y)=\sum_{\lambda_m\neq 0}\frac{\phi_m(x)\overline{\phi_m(y)}}{\lambda_m}$
and making use of (\ref{Bkernel}), one arrives at  (\ref{comp2}). $\blacksquare$

\begin{cor} For flat metrics $|\omega|^2$ with trivial holonomy the Green function $G_{hol}$ of the holomorphic extension $\Delta_{hol}$ is related to the Friedrichs Green function $G_F$ via
$$G_{hol}(x, y)=\int_X\partial_x G_F(x, z)\partial_{\bar y} G_F(z, y)\frac{1}{\omega(x)\overline{ \omega(y)}}d S(z)\,$$
\end{cor}

\subsection{Kernels of $\Delta_{hol}$ and $\Delta_{sing}$.}

\subsubsection{Kernel of the holomorphic extension} The following Proposition gives the complete description of the kernel of the holomorphic extension of the symmetric Laplacian on $X\setminus\{P\}$.
\begin{prop}\label{Weier} Suppose $P$ is not a Weierstrass point of $X$. Then the kernel of $\Delta_{hol}$ consists of constants: ${\rm dim}\,{\rm Ker}\, \Delta_{hol}=1$.
If $P$ is a Weierstrass point then the kernel of $\Delta_{hol}$ has dimension $2$ and is spanned by $1$ and a meromorphic function with  single pole at $P$ of multiplicity 2.
\end{prop}
{\bf Proof.} Let $U\in {\rm Ker}\,\Delta_{hol}$. Let $\xi$ be distinguished parameter near $P$ and let $X_{\epsilon}=X\setminus\{|\xi|\leq \epsilon\}$.  Using Stokes formula, one gets
$$0=-\frac{1}{4}<u, \Delta u>=\lim_{\epsilon\to 0}\Big\{\iint_{X_\epsilon}|\bar \partial u|^2+\oint_{|\xi|=\epsilon}\overline{(A/\xi^2+B/\xi+C+D\xi+E\xi^2+o(|\xi|^2))}\times$$
$$\partial_{\bar \xi}(A/\xi^2+B/\xi+C+D\xi+E\xi^2+o(|\xi|^2))\Big\}=\iint_X|\bar \partial u|^2$$
and, therefore,
$U$ is meromorphic on $X$ with a single pole of degree less or equal to $2$ at $P$.
It remains to notice that
\begin{itemize}
\item there are no meromorphic functions with a single pole of order $1$ on Riemann surfaces of positive genus,
\item for Riemann surfaces $X$ of genus 2 the point $P\in X$ is a Weierstrass point if and only if there exists a meromorphic function on $X$ with single double pole at $P$.
\end{itemize} $\square$

\subsubsection{Singular extension: very symmetric case}
Consider a hyperelliptic surface $X$ of genus 2 via $\mu^2=\prod_{j=1}^6(\lambda-\lambda_j)$ with  $\lambda_{k}=\lambda_1+r^2e^{\frac{2\pi i (k-1)}{5}}$;  $k=2, \dots, 6$; $r>0$.
Consider a holomorphic $1$-form $\omega$ on $X$ given by
$$\omega=\frac{(\lambda-\lambda_1)\,d\lambda}{\sqrt{\prod_{j=1}^6(\lambda-\lambda_j)}}\,.$$
Clearly, $\omega$ has a double zero at $P=(\lambda_1, 0)\in X$ and the metric $|\omega|^2$ is a flat metric on $X$ with unique conical point at $P$ of angle $6\pi$.
\begin{prop} The kernel of the singular self-adjoint extension $\Delta_{sing}$ of a symmetric Laplacian on $X\setminus \{P\}$ has dimension $3$.
\end{prop}
{\bf Proof.} There are two natural holomorphic local parameters on $X$ near $P$:
the one related to the ramified double covering $X\ni (\lambda, \mu)\mapsto \lambda \in {\mathbb C}\subset {\mathbb P}^1$,
$$\zeta=\sqrt{\lambda-\lambda_1}\,,$$
and the distinguished parameter $\xi$ for the conical metric $|\omega|^2$ related to the parameter $\zeta$ via
$$\xi(\zeta)=\left(\int_0^\zeta\frac{2w^2dw}{\sqrt{w^{10}-r^{10}}}\right)^{1/3}\,.$$
Since
$$\omega=C(\zeta^2 +O\left(\zeta^{12})\right)d\zeta\,$$
and, therefore,
\begin{equation}\label{sootn}\xi^3=C\left(\zeta^3+O(\zeta^{13})\right)\,,\end{equation}
one has
\begin{equation}\label{kvadrat}\frac{1}{\zeta^2}=\frac{C}{\xi^2}+o(|\xi|^3)\end{equation}
as $\xi, \zeta \to 0$ ($C$ is a constant which differs from one formula to another).
Now (\ref{kvadrat}) implies that
the meromorphic function
\begin{equation}\label{funnakr}P\mapsto f(P)=\frac{1}{\lambda(P)-\lambda_1}\end{equation}
on $X$ with a single double pole at $P$ belongs to ${\rm Ker}\Delta_{sing}$.
Clearly, $\bar f$ and $1$ also belong to ${\rm Ker}\Delta_{sing}$. Thus, ${\rm dim}\,{\rm Ker}\Delta_{sing}\geq 3$.

It turns out that in the case of the surface $X$ one can further specify the asymptotical expansion of the (unique up to a constant) harmonic function $g$ on $X$ with a single singularity at $P$ with
\begin{equation}\label{minusodin}g(\xi, \bar \xi)=\frac{1}{\xi}+O(1)\,.\end{equation}
Namely,
one has
\begin{equation}\label{nakonec}g=\frac{1}{\xi}+C+\alpha \bar \xi +o(|\xi|^2)   \end{equation}
with $\alpha\neq 0$.

According to
(\ref{hh2}) the coefficient near $\xi$
in the asymptotical expansion of $g$ near $P$ is equal to
$-\frac{1}{6}S_{Sch}(\xi)|_{\xi=0}$.
Using ${\mathbb Z}_5$-symmetry of $X$, it is easy to show that this coefficient must vanish.
First, notice that this quantity vanishes if
\begin{equation}\label{perehod1}S_{Sch}(\zeta)|_{\zeta=0}=0\end{equation}
This follows from the change of variables rule for
 a projective connection:
\begin{equation}\label{rule}S_{Sch}(\xi)=S_{Sch}(\zeta)\left(\frac{d\zeta}{d\xi}\right)^2+\{\zeta, \xi\}\end{equation}
(due to (\ref{sootn}) the Schwarzian derivative in the right hand side of the last equality vanishes at $\xi=0$).

Without loss of generality one can assume that $\lambda_1=0$. Consider the automorphism of $X$
$$\lambda\mapsto e^{\frac{2\pi i}{5}} \lambda\,.$$
Under this automorphism $\zeta\mapsto e^{\frac{\pi i}{5}}\zeta$ and, since the Schiffer projective connection is independent of the choice of basic cycles on $X$, one gets from (\ref{rule}) the relation
$$S_{Sch}(\zeta)|_{\zeta=0}=e^{\frac{2\pi i}{5}}S_{Sch}(\zeta)|_{\zeta=0}$$
implying $S_{Sch}(\zeta)|_{\zeta=0}=0$.

Using (\ref{TT}) and (\ref{funnakr}) together, it is easy to show that in the asymptotic expansion of $g$ there are no $\xi^2$ and $\bar \xi^2$ terms.
Notice that normalized holomorphic differentials $v_1$ and $v_2$ on $X$ are linear combinations of
$$\omega_1=\frac{d\lambda}{\sqrt{\prod_{k=1}^6(\lambda-\lambda_k)}}=\frac{2d\zeta}{\sqrt{\zeta^{10}-r^{10}}}$$
and
$$\omega_2=\frac{\lambda d\lambda}{\sqrt{\prod_{k=1}^6(\lambda-\lambda_k)}}=\frac{2\zeta^2 d\zeta}{\sqrt{\zeta^{10}-r^{10}}}$$
and, therefore,
$$v'_{1, 2}(\zeta)|_{\zeta=0}=0\,$$
Since $v'_{1, 2}(\xi)=v_{1, 2}'(\zeta)\left(\frac{d\zeta}{d\xi}\right)^2+v_{1, 2}(\zeta)\frac{d^2\zeta}{(d\xi)^2}$
and $\frac{d^2\zeta}{(d\xi)^2}=0$ (due to \ref{sootn})), one gets
\begin{equation}\label{diff0}v_{1, 2}'(\xi)|_{\xi=0}=0\,.\end{equation}
Relation (\ref{funnakr}) implies that one has $T_{12}(0)=0$ in (\ref{TT}) and from the symmetry $H(x, y)=H(y, x)$ of the function $H$ from Proposition \ref{szero} and (\ref{diff0}) one concludes that
$T_{21}(0)=0$ and,  therefore, there is no $\xi^2$ term in the expansion of $g$. Due to (\ref{diff0}) one has $T_{41}(0)=0$ and, therefore, the $\bar\xi^2$ is also absent.

It remains to notice that the coefficient $\alpha$ near $\bar\xi$ equals to $\pi B(\xi, \bar \xi)|_{\xi=0}$. Since  the imaginary part of the matrix of $b$-periods, $\Im {\mathbb B}$, is positively definite one has $\alpha\neq 0$ and (\ref{nakonec}) is proved.

To prove that ${\rm Ker}\,\Delta_{sing}={\rm lin.\ span}\{f, \bar f, 1\}$ it suffices to prove that a function $W$ from ${\rm Ker\,}\Delta^*$ with asymptotics
$$W=\frac{A}{\xi}+\frac{B}{\bar \xi}+o(|\xi|^2)$$
cannot belong to  ${\rm Ker}\,\Delta_{sing}$ unless $A=B=0$.
Assuming $W\in {\rm Ker}\,\Delta_{sing}$, one has $W-Ag-B\bar g\in {\rm Ker\,} \Delta_F$ and, therefore,
$$W=Ag+B\bar g+C$$
which contradicts (\ref{nakonec}) unless $A=B=0$. $\square$

\section{Appendix}

In this appendix, we sketch a proof of (\ref{MAINAS}). An alternative proof (of a closely related statement) based on different technical tools can be found in \cite{Moors}.



{\bf Step 1. Weighted Sobolev estimate for the infinite cone.}
Here we closely follow \cite{NazPlam}, Chapter 2, where a similar estimate was established for the Neumann boundary value problem in an angle lying in Euclidean plane. The only modification is the use of (\ref{RESKER}) instead of the Green function (2.14) from  \cite{NazPlam}.

Let $K$ be a standard round cone with conical point of conical angle $2\pi B$, $B=b+1>0$. Let $r>0$, $\omega$ ($0\leq \omega\leq 2\pi B$) be the standard polar coordinates on $K$.
Introduce the weighted Sobolev space $H^s_\beta(K)$ with norm
$$||u; H^l_\gamma (K)||=\left(\sum_{|\alpha|<l}\int_Kr^{2(\gamma-l+|\alpha|)}|\partial_x^\alpha u(x)|^2dx\right)^{1/2}$$
which is equivalent to the norm
$$\left(\sum_{j+k=l}\int_0^\infty\int_0^{2\pi B}r^{2(\beta+j)-1}|\partial^j_r\partial_\omega^k u(r, \omega)|^2\,dr\,d\omega \right)^{1/2}\,$$
where $$\gamma=\beta+l-1\,.$$

Using the Mellin transform $$u(r, \omega)\mapsto \tilde {u}(\lambda, \omega)=\frac{1}{\sqrt{2\pi}}\int_0^\infty r^{-i\lambda-1}u(r, \omega)\,dr\,,$$
one passes from  the
problem
\begin{equation}\label{ProbLap}-\Delta u=f \end{equation} in $K$
to the problem with parameter
\begin{equation}\label{probcircle}(-\partial_\omega^2+\lambda^2)\tilde{u}=(r^2 f)\tilde{\ }:=\tilde{F}\end{equation}
on the circle $S_{2\pi B}^1$
In what follows we need an explicit expression for the Green function of the problem (\ref{probcircle}). In order to obtain this expression
we start with Green function (integral kernel of the inverse operator)
$$\Phi(|x-y|)=\frac{\pi}{\lambda}e^{-\lambda|x-y|}$$
of the operator $\lambda^2-(d/dx)^2$ on ${\mathbb R}^1$ (see, \cite{Taylor}, vol. 1, (5.30) on p. 220).
Then the  Green  function of (\ref{probcircle}) is given by
$$\sum_{n\in {\mathbb Z}}\Phi (|\omega_1-\omega_2+2\pi B n|)\,.$$
 Summing geometric sequences, one immediately gets the needed expression for the Green function of (\ref{probcircle}):
\begin{equation}\label{RESKER}\Gamma(\omega_1, \omega_2; \lambda)=\frac{\pi}{\lambda}\frac{e^{-\lambda |\omega_1-\omega_2|}+e^{\lambda |\omega_1-\omega_2|}e^{-2\pi B}}{1-e^{-2\pi B \lambda}}\,.\end{equation}
Notice that $\Gamma$ has double pole at $\lambda=0$ and simple poles at $\lambda=\pm \frac{k}{B}$, $k=1, 2, \dots$

Clearly, one has the estimate
$$\int_0^{2\pi B}\int_0^{2\pi B}|\Gamma|^2 d\omega_1 d\omega_2\leq C |\lambda|^{-4}$$
for $\lambda\in {\mathbb R}+i\beta$; where $\beta$ is any real number not equal to  $\pm \frac{k}{B}$, $k=0, 1, 2, \dots$, and therefore,
one has an a priori estimate for (\ref{probcircle}):
 \begin{equation}\label{apr} ||\tilde{u}||^2\leq C |\lambda|^{-4} ||\tilde{F}||^2 \end{equation}
 with $C$ independent of $\lambda$  on the line ${\mathbb R}+i\beta$; $\beta\neq \pm \frac{k}{B}$, $k=0, 1, 2, \dots$
Using $\partial_\omega^{k+2}\tilde{u}=-\lambda^2\partial_\omega^k\tilde{u}-\partial_\omega^k\tilde{F}$, one upgrades (\ref{apr}) to
\begin{equation}\label{apr1}
\sum_{j=0}^{l+2}|\lambda|^{2(l+2-j)}||\partial_\omega^j\tilde{u}; L_2(S^1_{2\pi B})||^2\leq C \sum_{j=0}^l|\lambda|^{2(l-j)}||\partial_\omega^j\tilde{F}; L_2(S^1_{2\pi B})||^2\,.
\end{equation}
 on any horizontal line in ${\mathbb C}$ not passing through the poles of the Green function (\ref{RESKER}) (i. e. $\pm \frac{ki}{B}$, $k=0, 1, \dots$).
Returning from the Mellin images to the originals and using  estimate (\ref{apr1}) together with Parseval equality, one derives the following proposition.
\begin{prop}\label{estsol}Let $\gamma-l-1\neq \pm \frac{k}{B}$, $k=0, 1, 2, \dots$ and $f\in H^l_\gamma(K)$. Then there exists the unique solution $u\in H^{l+2}_\gamma(K)$ of  problem (\ref{ProbLap}).
One has the estimate
$$||u; H^{l+2}_\gamma(K)||\leq C ||f; H^l_\gamma(K)||\,.$$
\end{prop}

{\bf Step 2. Asymptotics for functions from ${\cal D}(\bar \Delta)$.}  Let $u\in C_0^\infty(X\setminus \{P\})$ and let $\chi$ be a smooth cut off function with a support in a small vicinity of $P$.
One can assume that $K$ and $X$ coincide in the latter vicinity.
Choose a small positive $\delta$ such that $\delta-1\neq \frac{k}{B}$; $k=0, 1, 2, \dots$; it should be noted that if $B$ is not an integer  (i. e. the conical angle is not an integer multiple of $2\pi$) one can take $\delta=0$.    Using Proposition \ref{estsol}, one gets
$$||\chi u; H^2_\delta (K)||\leq C ||\Delta(\chi u); H^0_\delta(K)||\leq C ||\Delta(\chi u); L_2(X)||\leq $$
$$C(||\chi \Delta u; L_2(X)||+ ||(\nabla \chi) \nabla u; L_2(X)||+||(\Delta \chi) u; L_2(X)||\,.$$
Since the support of $\nabla \chi$ is separated from $P$, one can estimate the second term in the right hand side from the above using the standard elliptic estimate (see, e. g, \cite{Berezin}, Theorem 2.1 of Appendix 2), thus, arriving at
\begin{equation}\label{estDOM}
||\chi u; H^2_\delta (K)||\leq C (||\Delta u; L_2(X)||+||u; L_2 (X)||)\,.
\end{equation}
The latter inequality shows that the functions from ${\cal D}(\bar \Delta)$ in a vicinity of $P$ belong to $H^2_\delta$ with (any small) positive $\delta$.
 Now the standard Sobolev lemma gives
 $$\sup_{1/2\leq|x|\leq 1}|u(x)|^2\leq C\sum_{|\alpha|\leq 2}\int_{1/2\leq |x|\leq 1} r^{2(\delta-2+|\alpha|)}|\partial^\alpha_x u(x)|^2dx$$
 with a constant $C$ independent of $u\in H^2_\delta$.
 Thus,
 $$\sum_{|\alpha|\leq 2}\int_{\epsilon/2\leq |x|\leq \epsilon} r^{2(\delta-2+|\alpha|)}|\partial^\alpha_x u(x)|^2dx=$$
 $$\epsilon^{2(\delta-1)}\sum_{|\alpha|\leq 2}\int_{1/2\leq |x|\leq 1} r^{2(\delta-2+|\alpha|)}|\partial^\alpha_x u(\epsilon x)|^2dx\geq
 C^{-1}\epsilon^{2(\delta-1)}\sup_{\epsilon/2\leq |x| \leq \epsilon}|u(x)|^2$$
 and, therefore,
\begin{equation}\label{ocenka}
u=O(r^{1-\delta})
\end{equation}
for $u\in {\cal D}(\bar \Delta)$ near $P$. The latter estimate can be improved to $u=O(r)$ in case of conical angles not equal to an integer multiple of $2\pi$.

In particular, (\ref{ocenka}) implies that all the terms in the right hand side of (\ref{MAINAS}) (except the last one) do not belong to ${\cal D}(\bar \Delta)$.

{\bf Step 3. Domain of the adjoint operator.} Let $u\in {\cal D}(\Delta^*)$. Then (see \cite{ReedSimon}, \S X.1)
$$u=u_1+u_2+u_3$$
with $u_{1, 2}\in {\rm Ker}\, (\Delta^*\pm i)$ and $u_3\in {\cal D}(\bar \Delta)$. It suffices to find the asymptotics of  functions from ${\rm Ker}\, (\Delta^*+ i)$. (The kernel of $\Delta^*-i$ can be described in the same way.)
So let $(\Delta^*+i)v=0$. This guarantees that $v\in L_2(X)$ and $v\in C^\infty(X\setminus\{P\})$ (due to the standard theory of elliptic equations). To establish the asymptotic behavior of $v$ at $P$ one has to put $v$ in some weighted Sobolev class. The following Lemma (going back to
\cite{Kondratjev}, \S 5.1) is actually the first key point of the whole proof.
\begin{lem} For some $\epsilon>0$ one has the inequality
$$\int_{\{x\in X\ : \ {\rm dist}(x, P)\leq \epsilon\}}r^4|\nabla^2 v|^2+r^2|\nabla v|^2 dx\leq \infty\,.$$
\end{lem}
{\bf Proof.} Let $\chi$ be the same cut off function as above and let $v_1=\chi v$.
Then \begin{equation}\label{Helm}\Delta v_1= -i v_1+f\end{equation} in $K$, where $f\in C^\infty_0(K\setminus \{P\})$.
Let $\phi, \kappa\in C_0^\infty (K\setminus\{P\})$; $\phi\kappa=\kappa$, ${\rm supp}\,\kappa \subset \{x: 1/4<|x|<4\}$. The standard elliptic estimate gives
\begin{equation}\label{standell}
||\kappa v_1; H^2(K)||\leq c\left(||\phi\Delta v_1; L_2(K)||+||\phi v_1; L_2(K)|| \right)\,.
\end{equation}
Choose  a partition of unity $\{\kappa_j\}$ and functions $\phi_j\in C^\infty_0(K\setminus\{P\})$ such that
$${\rm supp} \kappa_j\subset \{x: 2^{j-1}\leq |x|\leq 2^{j+1}\}\,,$$

$${\rm supp} \phi_j\subset \{x: 2^{j-2}\leq |x|\leq 2^{j+2}\}\,,$$
$\kappa_j\phi_j=\kappa_j$ and
$$|d^\alpha\kappa_j|+|D^\alpha\phi_j|\leq C_\alpha2^{-j|\alpha|}\,.$$
Under the rescaling $x\mapsto 2^j x$ the estimate (\ref{standell}) turns into
\begin{equation}\label{stand2}\sum_{|\alpha|\leq 2}2^{j(2|\alpha|-4)}\int_K|D^\alpha_x(\kappa_j v_1)|^2\,dx\leq C\left(\int_K|\phi_j\Delta v_1|^2\,dx +2^{-4j}\int_K|\phi_j v_1|^2\,dx  \right)\end{equation}
Taking into account (\ref{Helm}), multiplying (\ref{stand2}) by $2^{4j}$ and summing through $j=0, -1,$ $-2,...$, one gets the Lemma.

Thus, the solution $v_1$ of
the problem
\begin{equation}\label{prob3} \Delta w+iw = f \end{equation}
on the cone $K$ with $f\in C^\infty_0(K\setminus \{P\})$
belongs to the weighted Sobolev class $H^2_2(K)$.
It is known (see, e. g., \cite{MazNazPlam}, Section 1.3.6 for a similar statement for a problem in an angle, the needed modifications in the case of cone $K$ are obvious)  that two solutions $v_1$ and $v_2$ of (\ref{prob3}) from different weighted Sobolev classes $H^2_{\beta_1}(K)$ and $H^2_{\beta_2}(K)$; $\beta_1<\beta_2$ are related via
\begin{equation}\label{asy3}\chi v_1=\chi \sum_k W_k+\chi v_2\,\end{equation}
 where $\chi$ is the usual cut off function and $W_k$ are those special solutions $W_k(r, \omega)=a_k(r)e^{i\frac{k}{B}\omega}$ of the homogeneous problem
 $$\Delta w+ iw=0$$
on $K$ that satisfy $\chi  W_k\in H^2_{\beta_1}$ and $\chi W_k\notin H^2_{\beta_2}$ (there is always a finite number of such $W_k$). In fact, functions $a_k(r)$ can be expressed through Bessel functions. More precisely,
from (\ref{asy3})
one gets the
representation
\begin{equation}\label{predvas}
\chi v_1=\chi \sum_{k=-[B]}^{-1}K_{|k/B|}(e^{\frac{3\pi i}{4}}r)(c_k e^{i\frac{\omega k}{B}} +c'_k e^{-i\frac{\omega k}{B}})+ \chi c_0 K_0 (e^{\frac{3\pi i}{4}}r)+
\end{equation}
$$\chi c_0' I_0(e^{\frac{3\pi i}{4}}r) +   \chi \sum_{k=1}^{[B]}I_{|k/B|}(e^{\frac{3\pi i}{4}}r)(c_k e^{i\frac{\omega k}{B}} +c'_k e^{-i\frac{\omega k}{B}}) + \chi v_2\,$$
with $v_2=O(r^{1+\epsilon})$ with $\epsilon>0$.
Using \ref{predvas} and the series for Bessel functions (\cite{Lebedev}: (5.7.1), (5.7.2) and (5.7.11)) (here is the second key point: the powers of $r$ in the neighbour terms of the series for Bessel functions differ by two!),
one arrives at the asymptotic representation of $v_1$ of the type (\ref{MAINAS}) with the remainder $R=O(r^{1+\epsilon})$ which is smooth outside the vertex $P$ (moreover, this asymptotic expansion can be differentiated).
Let $\psi=1-\chi$, then it is straightforward to show that the sequence $\Delta[\psi(nx) R(x)]$ is uniformly bounded in $L_2(X)$ as $n\to\infty$ and, therefore, $R\in {\cal D}((\Delta^*)^*)= {\cal D}(\bar \Delta)$. $\square$

\end{document}